\documentclass[reqno]{amsart}
\usepackage{amsmath}
\usepackage{amsthm}
\usepackage{amsfonts}
\usepackage{amssymb}
\usepackage{latexsym}
\usepackage{amscd}
\usepackage{stmaryrd}
\usepackage{color}
\usepackage{hyperref}
\usepackage{extarrows}

\theoremstyle{plain}
\newtheorem{thm}{Theorem}[section]
\newtheorem{cor}[thm]{Corollary}
\newtheorem{lem}[thm]{Lemma}

\numberwithin{equation}{section}

\newfont{\scyr}{wncyr10 scaled 550}

\allowdisplaybreaks
\begin{document}

\title{Order structure and topological properties of the set of multiple $t$-values}

\author{Ende Pan}

\address{College of Teacher Education, QuZhou University, No. 78 jiuhua Road,
 QuZhou,PR China}

\email{1710385@tongji.edu.cn}

\keywords{multiple $t$-values, order structure, derived set.}

\subjclass[2010]{11M32, 06A05}

\begin{abstract}
In this paper, we compute the iterated derived sets of the set of multiple $t$-values under the usual topology of $\mathbb{R}$. Our results imply that the set of multiple $t$-values, ordered by $\geq$, is a well-ordered set. We determine its type of order, which is $\omega^{2}$, where $\omega$ is the smallest infinite ordinal. There exists a unique bijection from the set of multiple $t$-values to $\mathbb{N}^{2}$, which reverses the orders. We provide some description of this bijection.
\end{abstract}

\maketitle


\section{Introduction and statement of main results}\label{Sec:Intro}
For a finite sequence $\mathbf{k}=(k_1,\ldots,k_d)$ of positive integers with $k_{1}>1$, the multiple zeta value $\zeta(\mathbf{k})$ is defined by the following infinite series
$$\zeta(\mathbf{k})=\zeta(k_1,\ldots,k_d)=\sum\limits_{m_1>\cdots>m_d>0}\frac{1}{m_1^{k_1}\cdots m_d^{k_d}}.$$
In recent years there are some different variations of the multiple zeta values. In this paper we consider the multiple $t$-values, which were first systematically studied by Hoffman in \cite{Hoffman}.
Here a  multiple $t$-value $t(\mathbf{k})$ is defined as
\begin{align}t(\mathbf{k})=t(k_1,\ldots,k_d)=\sum\limits_{m_1>\cdots>m_d>0}\frac{1}{(2m_1-1)^{k_1}\cdots(2m_d-1)^{k_d}}.\end{align}
Usually $d$ is called the depth, $k_1+\cdots+k_d$ is called the weight, and $\mathbf{k}$ is called an admissible multi-index. By convention we denote the empty index by $\emptyset$, and define $t(\emptyset)$ to be $1$.

Previous research on multiple $t$-values revealed that they remarkable parallel to, and contrast with, multiple zeta values. In \cite{J. Zhao}, Zhao obtained some type of sum formula for multiple t-values. More precisely for positive integers $d$, $n$ with $d<n+1$, let $T(2n,d)$ be the sum of all multiple $t$-values with even arguments whose weights are $2n$ and whose depths are $d$. Then it was proved in \cite{J. Zhao} that
$$T(2n,d)=\frac{(-1)^{n-d}\pi^{2n}}{4^{n}(2n)!}\sum_{l=0}^{n-d}{\binom {n-l} d}{\binom {2n} {2l}}E_{2l},$$ where $E_{j}$ is Euler's number. More general weighted sums of multiple $t$-values at even arguments with polynomial weights were studied by Li and Xu in \cite{LI XU}. In \cite{Hoffman} Hoffman gave explicit formulas for the multiple $t$-values with repeated arguments analogous to those known for multiple zeta values by quasi shuffle products. For example, the well-known identities
  $$ \zeta({\{2\}}_{n})=\frac{\pi^{2n}}{(2n+1)!},\quad \zeta({\{4\}}_{n})=\frac{2^{2n+1}\pi^{4n}}{(4n+2)!}, \quad \zeta({\{6\}}_{n})=\frac{6{(2\pi)}^{6n}}{(6n+3)!}$$
  (where ${\{k\}}_{n}$ means $k$ repeats $n$ times) have multiple t-values counterparts
\begin{align}
t({\{2\}}_{n})=\frac{\pi^{2n}}{(2n)!2^{2n}},\quad t({\{4\}}_{n})=\frac{\pi^{4n}}{(4n)!2^{2n}},\quad t({\{6\}}_{n})=\frac{3{(\pi)}^{6n}}{(6n)!4}
. \end{align} Furthermore, Hoffman provided some interesting conjectures. For example he conjectured that the dimension of the rational vector space generated by weight $n$ multiple $t$-values equals to the $n$th Fibonacci number.

  For multiple zeta values, Kumar has a unconventional idea. In \cite{Senthil Kumar}, he studied the order structure and the topological properties of the set $\mathcal{Z}$ of all multiple zeta values. Taking the usual order and the usual topology of the set $\mathbb{R}$ of real numbers, Kumar computed the derived sets of the topological subspace $\mathcal{Z}$ of $\mathbb{R}$, and showed that the set $\mathcal{Z}$, ordered by $\geqslant$, is well-ordered with the order type $\omega^3$, where $\omega$ is the smallest infinite ordinal. Inspired by the work of Kumar, Li and the author studied the topological properties of some $q$-analogues of multiple zeta values in \cite{LI PAN}.

  In this paper, we study the order structure and the topological properties of multiple $t$-values. Let $\mathcal{T}$ be the set of all multiple $t$-values. Under the usual topology, we determine the sequence $(\mathcal{T}^{(n)})_{n\geqslant 0}$ of the derived sets of the topological subspace $\mathcal{T}$ of $\mathbb{R}$. Here $\mathcal{T}^{(0)}=\mathcal{T}$ and for any $n\in\mathbb{N}$, $\mathcal{T}^{(n)}$ is the set of accumulation points of $\mathcal{T}^{(n-1)}.$ As usually, we denote $\mathcal{T}^{(1)}$ by $\mathcal{T}'$ and $\mathcal{T}^{(2)}$ by $\mathcal{T}''.$
  To state our result about $\mathcal{T}^{(n)}$, we introduce the concept of the tail of multiple $t$-values. For an admissible multi-index $\mathbf{k}=(k_1,\ldots,k_d)$ and a positive odd integer $n$, we define
  $$t(\mathbf{k})_{n}=t(k_1,\ldots,k_d)_{n}=\sum\limits_{m_1>\cdots>m_d>n}\frac{1}{(2m_1-1)^{k_1}\cdots(2m_d-1)^{k_d}}.$$
  We set $t(\mathbf{\emptyset})_{n}=1.$ Then under the usual topology of $\mathbb{R}$, the derived sets of $\mathcal{T}$ are described as in the following theorem.
\begin{thm}\label{Thm:Derived-Set}
We have $$\mathcal{T}'=\{t(\mathbf{k})_{1}\mid \mathbf{k}\text{\;is admissible}\}\cup \{0,1\}$$ and $\mathcal{T}^{(n)}=\{0\}$ for any positive integer $n\geq 2$.
\end{thm}For the order structure of the set $\mathcal{T}$, we have the following two results.
\begin{thm}\label{Thm:well-ordered set}
The set $\mathcal{T}$, ordered by $\geq,$ is a well-ordered set with $t(2)$ as the maximum element.
\end{thm}

\begin{thm}\label{Thm:order structure}
The type of the order of $(\mathcal{T},\geq)$ is $\omega^{2}$, where $\omega$ is the smallest infinite ordinal.
\end{thm}

The paper is organized as follows. In Section 2, we study the derived sets of $\mathcal{T}$ and give a proof of Theorem \ref{Thm:Derived-Set}. In Section 3, we recall the definition of well-ordered sets and give a proof of Theorem \ref{Thm:well-ordered set}. In Section 4, we recall some properties of ordinals, and then give a proof of Theorem \ref{Thm:order structure}. By Theorem \ref{Thm:order structure}, there is a unique bijection $\Phi: (\mathcal{T},\geq)\rightarrow (\mathbb{N}^{2},\leq),$ which reverses the orders. We give some description of this bijection at the end of Section 4.

\section{Proof of  Theorem \ref{Thm:Derived-Set}}\label{Sec:Proof-Thm1-1}
In this section, we give a proof of  Theorem \ref{Thm:Derived-Set}.

We recall a formula.
\begin{lem}[\cite{Hoffman}]\label{Lem:identy}
Let $G=\sum_{j=0}^{\infty} \frac {(-1)^{j}}{2j+1}$ be Catalan's constant, then we have $$\sum_{j=1}^{\infty}t(2,\{1\}_{j-1})=2G.$$
\end{lem}
To prove Theorem \ref{Thm:Derived-Set}, we have to know the behaviour of the convergent sequences in the space $\mathcal{T}$. Recall that a sequence $(a(n))_{n\in\mathbb{N}}$ is injective if for any $n\neq m,$ it holds $a(n)\neq a(m)$.
\begin{thm}\label{Thm:accumulation point}
 Let $(\mathbf{k}(n))_{n\in\mathbb{N}}=((k_1(n),\ldots,k_{d(n)}(n)))_{n\in\mathbb{N}}$ be an injective sequence of admissible multi-indices. If $0$ is not the accumulation point of $(t(\mathbf{k}(n)))_{n\in\mathbb{N}}$, then there exist an infinite subset $\mathbb{D}$ of $\mathbb{N}$ and a positive integer $d$ such that
  \item [(i)] $d(n)=d$ for all $n\in\mathbb{D}$,
  \item [(ii)] and if $d>1,$ there are positive integers $k_1,\ldots,k_{d-1},$ with $k_i(n)=k_{i}$ for all $n\in\mathbb{D}$ and $1\leq i<d$.
\end{thm}
\noindent {\bf Proof.}
If the sequence $({d(n)})_{n\in\mathbb{N}}$ is unbounded, there exists a subsequence $({n_{k}})_{k\in\mathbb{N}}$ of $\mathbb{N}$ such that $\lim_{k\rightarrow\infty}d(n_{k})=\infty$. Then Lemma \ref{Lem:identy} implies that
$$\lim_{k\rightarrow\infty}t(2,\{1\}_{d(n_{k})-1})=0.$$
From the fact $$0<t(\mathbf{k}(n_{k}))\leq t(2,\{1\}_{d(n_{k})-1}),$$
we find that $0$ is an accumulation point of $(t(\mathbf{k}(n)))_{n\in\mathbb{N}},$ a contradiction. Hence there exist
an infinite subset $D_{1}$ of $\mathbb{N}$ and a positive integer $d$, such that $d(n)=d$ for all $n\in D_{1}$.
If $d=1$, we take $\mathbb{D}=D_{1},$ and complete the proof of the theorem.

Assume that $d>1$. If $({{k_{1}(n)}})_{n\in D_{1}}$ is unbounded, there exists a subsequence $n_{k}$ of $D_{1}$ with the property $$\lim_{k\rightarrow\infty}k_{1}(n_{k})=\infty.$$
Since $$t(\mathbf{k}(n_{k}))\leq t(k_{1}(n_{k}),\{1\}_{d-1}),$$ and
$$\lim\limits_{{k} \rightarrow \infty}t(k_{1}(n_{k}),\{1\}_{d-1})=0,$$
a contradiction. Hence there exist an infinite subset $D_{2}$ of $D_{1}$ and a positive integer $k_{1}$ such that $k_{1}(n)=k_{1}$ for all
$n\in D_{2}$. Now assume that for $1\leq j< d-1$, we have found an infinite subset $D_{j}$ of $\mathbb{N}$ and positive integers $k_1,\ldots k_{j}$, such that $k_{i}(n)=k_{i}$ for  $i=1,2,\ldots,j$ and all $n\in D_{j}$. Then using the inequality
$$t(k_1,\ldots k_{j},k_{j+1}(n),\ldots, k_d(n))<t(k_1,\ldots, k_{j})t(k_{j+1}(n),\ldots k_d(n))$$
with $k_{j+1}(n)>1,$ we conclude that there are an infinite subset $D_{j+1}$ of $D_{j}$ and a positive integer $k_{j+1}$ such that $k_{j+1}(n)=k_{j+1}$ for all $n\in D_{j+1}.$ Finally we take  $\mathbb{D}=D_{d-1}$ and finish the proof. \qed

We can restate Theorem \ref{Thm:accumulation point} in the following way.
\begin{cor}\label{Cor:injective point-3}Let $(\mathbf{k}(n))_{n\in\mathbb{N}}$ be an injective sequence of admissible multi-indices such that $0$ is
not an accumulation point of $(t(\mathbf{k}(n)))_{n\in\mathbb{N}}$. Then $(\mathbf{k}(n))_{n\in\mathbb{N}}$ has a subsequence of the following type
\begin{align}({\mathbf{k},{\varphi(n)+2}})_{n\in \mathbb{N}},\end{align}\label{Eq:Type-2}where $\mathbf{k}$ is a fixed admissible multi-index or an empty index and $({{\varphi(n)}})_{n \in \mathbb{N}}$ is a strictly increasing sequence in $\mathbb{N}$.
\end{cor}

Therefore we have the following corollaries.
\begin{cor}\label{Cor:injective point-4}
Each injective sequence $(\mathbf{k}(n))_{n\in\mathbb{N}}$ of admissible multi-indices has a
subsequence $(\mathbf{k}(n_{k}))_{k\in\mathbb{N}}$ such that the sequence $(t(\mathbf{k}(n_{k})))_{k\in\mathbb{N}}$ is strictly decreasing.
\end{cor}
\noindent {\bf Proof.}
If $0$ is an accumulation point of $(t(\mathbf{k}(n)))_{n\in\mathbb{N}}$, we obviously get the result. Otherwise, the result follows from Corollary \ref{Cor:injective point-3}. \qed

\begin{cor}\label{Cor:injective point-1}
For any real number $\alpha,$ there exist only finitely many admissible multi-indices $\mathbf{k}$ for which
$t(\mathbf{k})=\alpha.$
\end{cor}
\noindent {\bf Proof.}
We immediately get the result from Corollary \ref{Cor:injective point-4}. \qed

Now we come to prove Theorem \ref{Thm:Derived-Set}

\noindent{\bf Proof of Theorem \ref{Thm:Derived-Set}.} We first compute $\mathcal{T}'$. For an admissible multi-index $\mathbf{k}=(k_{1},\ldots,k_{d}),$
we take $\mathbf{k}(n)=(k_{1},\ldots, k_d,n+2).$ Then since $$\lim_{n\rightarrow\infty}t(\mathbf{k}(n))=t(k_{1},\ldots, k_{d})_{1},$$ we get $t(\mathbf{k})_{1}\in \mathcal{T}'$.
Similarly, taking $\mathbf{k}(n)=(2,\{1\}_{n-1})$ and $\mathbf{k}(n)=(n+1)$ for $n\in \mathbb{N}$ respectively, we find $0$ and $1$ belong to $\mathcal{T}'$.

Conversely, for any nonzero $\alpha\in \mathcal{T}'$, there is a sequence $(\mathbf{k}(n))_{n\in\mathbb{N}}$
such that $\lim_{n\rightarrow\infty}t(\mathbf{k}(n))=\alpha.$ By Corollary \ref{Cor:injective point-3} without loss of generality we may assume that $(\mathbf{k}(n))_{n\in\mathbb{N}}$ is one of the following types
 \begin{enumerate}
                \item  $({{\varphi(n)+2}})_{n\in \mathbb{N}}$;
                \item  $({\mathbf{k},{\varphi(n)+2}})_{n\in \mathbb{N}},$
              \end{enumerate}
where $\mathbf{k}=(k_{1},\ldots, k_d)$ is a fixed admissible multi-index and $({{\varphi(n)}})_{n \in \mathbb{N}}$ is a strictly increasing sequence in $\mathbb{N}$. In the case of type (1), since
\begin{align*}&\lim_{n\rightarrow\infty}t(\varphi(n)+2)\\
 &=\lim_{n\rightarrow\infty}\sum\limits_{m>0}\frac{1}{(2m-1)^{\varphi(n)+2}}\\
 &=\sum\limits_{m>0}\lim_{n\rightarrow\infty}\frac{1}{(2m-1)^{\varphi(n)+2}}\\
 &=1,\end{align*}
we get $\alpha=1$. And in the case of type (2), since
\begin{align*}&\lim_{n\rightarrow\infty}t(k_{1},\ldots, k_d,\varphi(n)+2)\\
&=\lim_{n\rightarrow\infty}\sum\limits_{m_1>\cdots>m_d>m_{d+1}>0}\frac{1}{(2m_1-1)^{k_1}\cdots(2m_d-1)^{k_d}(2m_{d+1}-1)^{\varphi(n)+2}}\\
&=\sum\limits_{m_1>\cdots>m_d>m_{d+1}>0}\lim_{n\rightarrow\infty}\frac{1}{(2m_1-1)^{k_1}\cdots(2m_d-1)^{k_d}(2m_{d+1}-1)^{\varphi(n)+2}}\\
&=\sum\limits_{m_1>\cdots>m_d>1}\frac{1}{(2m_1-1)^{k_1}\cdots (2m_d-1)^{k_d}}\\
&=t(k_{1},\ldots, k_{d})_{1},\end{align*}
we have $\alpha=t(k_{1},\ldots, k_{d})_{1}.$ Hence the result about $\mathcal{T}'$ is proved.

Now we compute $\mathcal{T}''$. Note that for any admissible multi-index $\mathbf{k}$, we have $t(\mathbf{k})> t(\mathbf{k})_{1}$. Hence we get $0\in \mathcal{T}''$ from $0\in \mathcal{T}'$.
Conversely, assume that $0$ is not the accumulation point of $(t(\mathbf{k}(n))_{1})_{n\in\mathbb{N}}.$
Then $0$ is also not the accumulation point of $(t(\mathbf{k}(n)))_{n\in\mathbb{N}}.$ Hence from Corollary \ref{Cor:injective point-3}, we may assume that $(\mathbf{k}(n))_{n\in\mathbb{N}}$ is of type (1) or of type (2).
Then from the facts
\begin{align*}&\lim_{n\rightarrow\infty}t(\varphi(n)+2)_{1}\\
&=\lim_{n\rightarrow\infty}\sum\limits_{m>1}\frac{1}{(2m-1)^{\varphi(n)+2}}\\
 &=\sum\limits_{m>1}\lim_{n\rightarrow\infty}\frac{1}{(2m-1)^{\varphi(n)+2}}\\
 &=0\end{align*} and
 \begin{align*}&\lim_{n\rightarrow\infty}t(k_{1},\ldots, k_d,\varphi(n)+2)_{1}\\
&=\lim_{n\rightarrow\infty}\sum\limits_{m_1>\cdots>m_d>m_{d+1}>1}\frac{1}{(2m_1-1)^{k_1}\cdots (2m_d-1)^{k_d}(2m_{d+1}-1)^{\varphi(n)+2}}\\
&=\sum\limits_{m_1>\cdots>m_d>m_{d+1}>1}\lim_{n\rightarrow\infty}\frac{1}{(2m_1-1)^{k_1}\cdots (2m_d-1)^{k_d}(2m_{d+1}-1)^{\varphi(n)+2}}\\
&=0,\end{align*}
we have $\mathcal{T}''=\{0\}.$

Finally, it is obvious that $\mathcal{T}^{(n)}=\{0\}$ for any $n\geq 2$. \qed

Similar as Corollary \ref{Cor:injective point-1}, we have the following result.
\begin{cor}\label{Cor:injective point-5}
For any real number $\beta,$ there exist only finitely many admissible multi-indices $\mathbf{k}$ for which
$t(\mathbf{k})_{1}=\beta.$
\end{cor}


\section{Proof of Theorem \ref{Thm:well-ordered set} }\label{Sec:Proof-Thm1-2}
In this section, we give a proof of  Theorem \ref{Thm:well-ordered set}.
We recall the definition of a well-ordered set. A partial ordered set $(X,\geq)$ is called well-ordered, if the order $\geq$ is totally ordered, and any nonempty subset of $X$ has a maximal element.

\noindent {\bf Proof of Theorem \ref{Thm:well-ordered set}.}
As a subset of $\mathbb{R}$, $\mathcal{T}$ is obviously totally ordered.
Now we prove that any nonempty subset $X$ of $\mathcal{T}$ has a maximal element.
If $X$ is a finite set, it is easy to find the maximal element of $X$ as $X$ is totally ordered. Now we investigate another case, $X$ is an infinite set. Assume that there is no maximal element of $X$. Then for any $\alpha_{1}\in X$,
there is a $\alpha_{2}\in X$, such that $\alpha_{1}<\alpha_{2}.$ Similarly there exists $\alpha_{3}\in X$ with $\alpha_{2}<\alpha_{3}$. Hence there is an infinite strictly increasing sequence $\alpha_{1}<\alpha_{2}<\alpha_{3}<\cdots$ in $X$, which contradicts Corollary \ref{Cor:injective point-4}.

We next prove that $t(2)$ is the maximal element of $\mathcal{T}$.
Since $t(2)=\frac{\pi^{2}}{8}$, from Lemma \ref{Lem:identy}, we obtain
$$\sum_{j=1}^{\infty}t(2,\{1\}_{j})=2G-\frac{\pi^{2}}{8}<1.$$ Hence we find $$t(2)>t(2,\{1\}_{j})$$ for all $j\in \mathbb{N}.$ Finally with the inequality $t(k_{1},\ldots,k_{j+1})\leq t(2,\{1\}_{j})$, we complete our proof.
\qed
\begin{cor}\label{Cor:well-ordered}
$\mathcal{T}'$ is well-ordered with $t(\emptyset)_{1}$ as the maximum element and $t(2)_{1}$ as the second largest element.
\end{cor}
\noindent {\bf Proof.}
Similar as the proof of Theorem \ref{Thm:well-ordered set}, to show that $\mathcal{T'}$ is well-ordered, we need the fact that for any injective sequence $(\mathbf{k}(n))_{n\in\mathbb{N}}$ of admissible multi-indices, there is a strictly decreasing infinite subsequence of $({t(\mathbf{k}(n))_{1}})_{n\in\mathbb{N}}$. While one can get this fact similar as Corollary \ref{Cor:injective point-4}.

We next compute the maximum element and the second largest element of $\mathcal{T'}$. From the equations
$$t(2)_{1}=t(2)-1$$ and
$$t(2,\{1\}_{j})_{1}=t(2,\{1\}_{j})-t(2,\{1\}_{j-1})_{1}\quad(j\geq 1),$$ we obtain
 $$\sum_{j=1}^{\infty}t(2,\{1\}_{j-1})_{1}=\sum_{j=1}^{\infty}t(2,\{1\}_{j-1})-1-\sum_{j=1}^{\infty}t(2,\{1\}_{j-1})_{1}.$$
So
$$\sum_{j=1}^{\infty}t(2,\{1\}_{j-1})_{1}=\frac{2G-1}{2}\approx0.416,$$ Hence we have $$t(\emptyset)_{1}=1>t(2,\{1\}_{j-1})_{1}\geq t(k_{1},\ldots,k_{j})_{1}$$ for all $j\in \mathbb{N}$
and we find $t(\emptyset)_{1}$ is the maximum element. Similarly, since $$t(2)_{1}\approx 0.232>\sum_{j=1}^{\infty}t(2,\{1\}_{j})_{1}>t(2,\{1\}_{j})_{1}\geq t(k_{1},\ldots,k_{j+1})_{1}$$ for all $j\in \mathbb{N},$ we have $t(2)_{1}$ is the second largest element. \qed

\section{Proof of  Theorem \ref{Thm:order structure}}\label{Sec:Proof-Thm1-3}
In this section, we first give a proof of Theorem \ref{Thm:order structure}.
We recall some lemmas.
\begin{lem}[\cite{Senthil Kumar}]\label{Lem:countable}
 Any well-ordered subset of  $\overline{\mathbb{R}}$ is countable, where $\overline{\mathbb{R}}$ is the extended real line.
\end{lem}
\begin{lem}[\cite{Senthil Kumar}]\label{Lem:ordinal}
Let $\mathcal{A}$ be a well-ordered subset of $\overline{\mathbb{R}}$. If $\mathcal{A}$ has order type $\omega\mu+\nu,$ where $\mu$ is an ordinal and $\nu$ is a finite ordinal, then the order type of $\mathcal{A}'$ is $\mu$ if $\mu$ is finite and $\mu + 1$ if $\mu$ is infinite.
\end{lem}

\noindent {\bf Proof of Theorem\ref{Thm:order structure}.}
Recall the division algorithm of ordinals. For ordinals $\gamma,\alpha,\beta$ with $\gamma<\alpha\beta$,
there exist unique ordinals $\alpha',\beta'$, such that $\gamma=\alpha\beta'+\alpha',$ and $\alpha'<\alpha$, $\beta'<\beta.$ Since $\mathcal{T'}$ is countable, we may assume the ordinal of $\mathcal{T'}$ is $\omega\mu+\nu,$ where $\mu$ is an ordinal and $\nu$ is a finite ordinal. From Lemma \ref{Lem:ordinal}, the ordinal of
$\mathcal{T''}$ is $\mu$ or $\mu+1$. But $\mathcal{T''}=\{0\}$, which implies that the ordinal of $\mathcal{T''}$ is $1$. Hence we get $\mu=1$. Since $0$ is the smallest element of $\mathcal{T'},$
we find $\nu=1$. Now we assume the ordinal of $\mathcal{T}$ is $\omega\mu'+\nu',$
where $\mu'$ is an ordinal and $\nu'$ is a finite ordinal. From Lemma \ref{Lem:ordinal} and the fact that the ordinal of $\mathcal{T'}$ is $\omega+1,$ we have $\mu'=\omega$. Since there is no smallest element
in $\mathcal{T}$, we get $\nu'=0$. Therefore, the ordinal of $\mathcal{T}$ is $\omega^{2}.$ \qed

Theorem \ref{Thm:order structure} guarantees that there exists a unique bijection
$$\Phi:\quad (\mathcal{T},\geq)\rightarrow (\mathbb{N}^{2},\leq),$$
reverses the orders. Here $\mathbb{N}^{2}$ is endowed with the lexicographical order. At the end of this section,
we provide some description of the map $\Phi$.

Since the ordinal of $\mathcal{T}'$ is $\omega+1$, there exists a unique bijection
$$\Psi:\quad (\mathcal{T'},\geq)\rightarrow (\mathbb{\overline{N}},\leq),$$
which reverses the orders. Here $\mathbb{\overline{N}}=\mathbb{N} \bigcup \{\infty\}$ and $\infty$ is the maximal element of $\mathbb{\overline{N}}$. Since $0$ is the minimal element of $\mathcal{T}'$, we have $\Psi(0)=\infty$.
And it is easy to compute the image of nonzero elements of $\mathcal{T}'$ under the map $\Psi$.
\begin{lem}\label{Lem:ordinal derived }
For any nonzero $\beta \in \mathcal{T}'$, we have $\Psi(\beta)=card(B_{\beta})$, where $$B_{\beta}=\{\alpha\in \mathcal{T}'| \alpha\geq \beta\}.$$
\end{lem}

Therefore, we may index the set $\mathcal{T}'-\{0\}$ by $\mathbb{N}$ as
$$\mathcal{T}'-\{0\}=\{\beta_{1}>\beta_{2}>\cdots>\beta_{n}>\cdots\}.$$
We also set $\beta_{0}=\infty.$
For any $r\in \mathbb{N}$, we set $$I_{r}=\{\mathbf{k}|\mathbf{k} \text{\; is admissible, \;} t(\mathbf{k})_{1}=\beta_{r}\}.$$
Which is a finite set by Corollary \ref{Cor:injective point-5}.
Set
$$\mathcal{L}_{r}=\{(\mathbf{k},n)|\mathbf{k}\in I_{r},n\in\mathbb{N},t(\mathbf{k},n)<\beta_{r-1}\}$$
and
$$\mathcal{P}_{r}=\{(\mathbf{k},n)|\mathbf{k}\in I_{r},n\in\mathbb{N},t(\mathbf{k},n)\geq \beta_{r-1}\}.$$
It is easy to show that $\mathcal{P}_{r}$ is finite for any $r\in \mathbb{N}$, and in particular $\mathcal{P}_{1}=\emptyset.$
For any $\alpha \in \mathcal{T}$, we set
$$\Phi(\alpha)=(\Phi_{L}(\alpha),\Phi_{R}(\alpha))\in\mathbb{N}^{2}.$$
Then we have the following lemma.
 \begin{lem}\label{Lem:iso-L }
For any $\mathbf{l}\in \mathcal{L}_{r},$ we have $\Phi_{L}(t(\mathbf{l}))=r.$
\end{lem}
 \noindent {\bf Proof.}
 We prove it by induction on $r$. For any $\mathbf{l}=(\mathbf{k},n)\in \mathcal{L}_{1}$,
 if $\Phi_{L}(t(\mathbf{l}))>1,$ then we have $$(1,1)<(1,2)<\cdots<(1,n)<\cdots<\Phi(t(\mathbf{l})).$$ Applying $\Phi^{-1}$, we get an injective sequence $(\mathbf{k}(n))_{n\in\mathbb{N}}$, such that $$t(\mathbf{k}(1))>t(\mathbf{k}(2))>\cdots>t(\mathbf{k}(n))>\cdots>t(\mathbf{l}).$$
 Therefore there exists an element $\beta\in \mathcal{T}'$ such that $$\beta\geq t(\mathbf{l}).$$
Hence $$\beta_{1}=t(\mathbf{k})_{1}<t(\mathbf{k},n)=t(\mathbf{l})\leq\beta,$$
which is impossible.

Now assume that $r>1$, and the result holds for any positive integer less than $r$.
Then for any $\mathbf{l}=(\mathbf{k},n)\in \mathcal{L}_{r-1}$, we have $\Phi(t(\mathbf{l}))=\Phi(t(\mathbf{k},n))=(r-1,\tau(n)).$
Since $$t(\mathbf{k},n)>t(\mathbf{k},n+1),$$
we have $({{\tau(n)}})_{n \in \mathbb{N}}$ is a strictly increasing sequence in $\mathbb{N}$.
For any $\mathbf{l'}\in \mathcal{L}_{r}$, we set $\Phi(t(\mathbf{l'}))=(r',\tau).$
We first prove $r'\geq r$. If $r'<r-1,$ then for any $\mathbf{l}\in \mathcal{L}_{r-1}$, we have
$$\Phi(t(\mathbf{l'}))<\Phi(t(\mathbf{l})),$$
which implies $$t(\mathbf{l'})>t(\mathbf{l})>\beta_{r-1},$$ a contradiction.
If $r'=r-1$, since $({{\tau(n)}})_{n \in \mathbb{N}}$ is a strictly increasing sequence in $\mathbb{N}$, we can find a big enough positive integer $m$ such that $\tau(m)>\tau$. Then there exists a multi-index $\mathbf{l}=(\mathbf{k},m)\in \mathcal{L}_{r-1}$ such that
$$\Phi(t(\mathbf{l'}))<\Phi(t(\mathbf{l})),$$
which implies $$t(\mathbf{l'})>t(\mathbf{l})>\beta_{r-1},$$ a contradiction.
We next prove $r'\leq r$. If $r'>r$, then $$(r,1)<(r,2)<\cdots<\Phi(t(\mathbf{l'})).$$
Applying $\Phi^{-1},$ there exists an injective sequence of admissible multi-indices $({\mathbf{k}(n)}_{n\in\mathbb{N}}$  such that $$t(\mathbf{k}(1))>t(\mathbf{k}(2))>\cdots>t(\mathbf{k}(n))>\cdots>t(\mathbf{l'}).$$ Since $$card(\bigcup_{1\leq i\leq r-1}\mathcal{P}_{i})$$ is finite, there exists an accumulation point $\beta$ of $(t(\mathbf{k}(n)))_{n\in\mathbb{N}}$ such that
 $$\beta \in \mathcal{T}'-\{\beta_{1},\ldots,\beta_{r-1}\}$$ and $\beta\geq t(\mathbf{l'})$, we get a contradiction. \qed

While for elements in $\mathcal{P}_{r}$, we have the following lemma.
 \begin{lem}\label{Lem:isolate point }
For any $\mathbf{l}\in \mathcal{P}_{r}$, we have $\Phi_{L}(t(\mathbf{l}))=\min\{m\in \mathbb{N}|t(\mathbf{l})>\beta_{m}\}.$
\end{lem}
\noindent {\bf Proof.}
We denote by $m$ the minimal number such that $\beta_{m}<t(\mathbf{l}).$ Then we have $\beta_{m-1}\geq t(\mathbf{l}).$
We first prove $\Phi_{L}(t(\mathbf{l}))\geq m$. If $m=1$, the conclusion is obvious. If $m>1$, by Lemma \ref{Lem:iso-L }
there exists an infinite sequence
$$({t(\mathbf{l}_{m-1})(n))})_{n\in\mathbb{N}}$$ such that for any positive integer $n$ $$t(\mathbf{l}_{m-1}(n))>t(\mathbf{l}_{m-1}(n+1))\geq t(\mathbf{l}),$$ where $\mathbf{l}_{m-1}\in \mathcal{L}_{m-1}.$
Therefore we obtain $$(m-1,\varsigma(1))<(m-1,\varsigma(2))<\ldots<(m-1,\varsigma(n))<\ldots<\Phi(t(\mathbf{l})),$$
where $({{\varsigma(n)}})_{n \in \mathbb{N}}$ is a strictly
increasing sequence in $\mathbb{N}$. Hence $\Phi_{L}(t(\mathbf{l}))>m-1$.

We next prove $\Phi_{L}(t(\mathbf{l}))<m+1$. If $\Phi_{L}(t(\mathbf{l}))\geq m+1$, then $$(m,1)<(m,2)<\cdots<\Phi(t(\mathbf{l})),$$
which implies that for any $\mathbf{l}_{m}\in \mathcal{L}_{m},$ it holds $\Phi(t(\mathbf{l}_{m}))>\Phi(t(\mathbf{l}))$. Then we have $\beta_{m}\geq t(\mathbf{l})$, which leads to a contradiction. \qed

From Lemma \ref{Lem:iso-L } and Lemma \ref{Lem:isolate point }, we find that for any $\alpha\in \mathcal{T}$ and $r\in \mathbb{N}$, we have $\Phi_{L}(\alpha)=r\Leftrightarrow\beta_{r-1}\geq\alpha>\beta_{r}.$ In other words we can order the subset $\mathcal{T}_{r}=\{\alpha\in \mathcal{T}|\beta_{r-1}\geq\alpha>\beta_{r}\}$ by $\mathbb{N}$ as $\{\alpha_{r,1}>\alpha_{r,2}>\cdots\}.$

Now we can get a description of the images of the map $\Phi$.
\begin{thm}\label{Thm:ordinal isomorphism}
For any $\alpha\in \mathcal{T}$, we have the following equality: $$\Phi(\alpha)=(\min\{i\in \mathbb{N}|\alpha>\beta_{i}\},\min\{j\in \mathbb{N}|\alpha\geq \alpha_{i,j}\}).$$
\end{thm}
\noindent {\bf Proof.}
We get this claim from the facts that there exists a unique bijection from $(\mathcal{T}_{i},\geq)$ to $(\mathbb{N},\leq)$ which reverses the orders, and there is a unique bijection $(\{(i,1),(i,2),\cdots\},\leq)$ to $(\mathbb{N},\leq)$ which preserves the orders. \qed

By the above description of the bijection
$\Phi:(\mathcal{T},\geq) \rightarrow (\mathbb{N}^{2},\leq)$
and the supplementary definition $t(1)=\infty$, we compute some approximate values and give the following ordering rule of $\mathcal{T}\bigcup\mathcal{T'}$:
$$t(1)>t(2)>t(3)>\cdots>\underline{t(\emptyset)_{1}}>
t(2,1)>t(2,2)>\cdots>\underline{t(2)_{1}}>
$$$$t(2,1,1)>t(2,1,2)>\cdots>\underline{t(2,1)_{1}}>
t(3,1)>t(3,2)>{\cdots>\underline{t(3)_{1}}}>\cdots.$$
Hence here we conjecture $\bigcup_{r=1}^{\infty}\mathcal{P}_{r}=\emptyset$ and there are no two different admissible multi-indices $\mathbf{k}$ and $\mathbf{l}$, such that $$t(\mathbf{l})_{1}=t(\mathbf{k})_{1}.$$
Moreover combining with Lemma \ref{Lem:iso-L } and Theorem \ref{Thm:ordinal isomorphism}, we have the following conjecture.

\noindent {\bf Conjecture.}
For any multi-index $\mathbf{l}=(\mathbf{k},n)$, we have $$\Phi(t(\mathbf{k},n))=(m,n),$$ where $\mathbf{k}$ is admissible or empty and $m=\min\{i\in \mathbb{N}|t(\mathbf{k})_{1}>\beta_{i}\}$.

\end{document}